\documentclass[10pt]{article}
\usepackage{amssymb,amsmath}
\usepackage{verbatim}
\usepackage[dvipsnames*,svgnames]{xcolor}
\usepackage{authblk}
\usepackage{graphicx}
\usepackage{subfigure}
\usepackage{listings}
\newcommand{\Keywords}[1]{\par\noindent
{\small{\em Keywords\/}: #1}}
\newtheorem{thm}{Theorem}
\newtheorem{lem}[thm]{Lemma}
\newtheorem{prop}[thm]{Proposition}
\newtheorem{dfn}[thm]{Definition}
\newtheorem{rem}[thm]{Remark}

\newenvironment{proof}[1][Proof]{\textbf{#1.} }{\ \rule{0.5em}{0.5em}}

\def\sgn{\mathop{\rm sgn}}
\let\eps=\varepsilon
\let\a=\alpha
\let\b=\beta
\let\D=\Delta
\let\d=\delta
\let\s=\sigma
\let\g=\gamma
\let\t=\theta

\let\kap=\varkappa
\let\L=\Lambda

\let\ls=\leqslant
\let\gs=\geqslant

\def\dvit{\colon\ }

\def\pr{\hbox{\bf P}}
\def\E{\hbox{\bf E}}

\date{}
\title{Exact confidence intervals of the extended Orey index for Gaussian processes}
\author[1,*,$\dag$]{K. Kubilius}
\author[2]{D. Melichov}
\affil[1]{Vilnius University, Institute of Mathematics and
Informatics, Akademijos 4, LT-08663, Vilnius, Lithuania}
\affil[2]{Vilnius Gediminas Technical University, Saul\.etekio al.
11, LT-10223, Vilnius, Lithuania}
\begin{document}

\maketitle

\let\oldthefootnote\thefootnote
\renewcommand{\thefootnote}{\fnsymbol{footnote}}
\footnotetext[1]{Corresponding author. E-mail:
\texttt{kestutis.kubilius@mii.vu.lt}\\ $\dag$ This research was funded by a grant (No.
MIP-048/2014) from the Research Council of Lithuania.}
\let\thefootnote\oldthefootnote


\abstract{In this paper exact confidence intervals for the Orey index of Gaussian processes are obtained using concentration inequalities for Gaussian quadratic forms and discrete observations of the underlying process. The obtained result is applied to Gaussian processes with the Orey index which not necessarily have stationary increments.

\bigskip
 }


\section{Introduction}

Let $X = \{X(t)\dvit t \in [0, T ]\}$ be a second order stochastic process with the
incremental variance function $\s_X^2$ defined on $[0, T ]^2 := [0, T ]\times[0, T ]$ with values
\[
\s_X^2(s, t) := \E[X(t)-X(s)]^2,\quad (s, t)\in[0, T ]^2.
\]

Denote by $\widetilde \Psi$ a class of continuous functions $\varphi\dvit (0,T]\to [0,\infty)$ such that $\lim_{h\downarrow  0}\varphi(h)=0$ and  $L(h)=\varphi(h)/h\to\infty$, $h\downarrow 0$. For example, we can take $\varphi(h)=h \vert \ln h\vert$ or $\varphi(h)=h^{1-\eps}$ for $0<\eps<1$.
Set
\begin{align*}
\g_*:=&\inf\bigg\{\g>0\dvit \lim_{h\downarrow 0}\sup_{\varphi(h)\ls s\ls T-h}\frac{h^\g}{\s_X(s,s+h)}=0\bigg\},\\
\widetilde\g_*:=&\inf\Big\{\g>0\dvit \lim_{h\downarrow 0}\frac{h^\g}{\s_X(0,h)}=0\Big\}
\end{align*}
and
\begin{align*}
\g^*:=&\sup\bigg\{\g>0\dvit \lim_{h\downarrow 0}\inf_{\varphi(h)\ls s\ls T-h}\frac{h^\g}{\s_X(s,s+h)}=+\infty\bigg\},\\
\widetilde\g^*:=&\sup\Big\{\g>0\dvit \lim_{h\downarrow 0}\frac{h^\g}{\s_X(0,h)}=+\infty\Big\}\,,
\end{align*}
where $\varphi\in\widetilde \Psi$. Note that $0\ls\widetilde\g^*\ls\widetilde\g_*\ls +\infty$ and $0\ls\g^*\ls\g_*\ls +\infty$. In paper \cite{kk1} we used a narrower class of functions $\Psi$, i.e. such functions which additionally satisfy condition $\lim_{h\downarrow  0}[h\cdot L^3(h)]=0$. This condition is not necessary for the existence of the Orey index. It is required only for consideration of almost sure asymptotic behavior of the second-order quadratic variations of Gaussian processes.

We give the following extension of the Orey index.
\begin{dfn}[\cite{kk1}]\label{oreydef}
Let $X = \{X(t)\dvit t \in [0, T ]\}$ be a second order stochastic process with the
incremental variance function $\s_X^2$ such that $\sup_{0\ls s\ls T-h}\s_X(s,s+h)\to 0$ as $h\to 0$. If $\g_*=\widetilde\g_*=\g^*=\widetilde\g^*$ for any function $\varphi\in\widetilde \Psi$, then we say that the process $X$ has the Orey index $\g_X=\g_*=\widetilde\g_*=\g^*=\widetilde\g^*$.
\end{dfn}

Assume that for some $\g\in(0,1)$ the second order stochastic process $X$ satisfies conditions:

(C1)\quad $\s_X(0,\d)\asymp\d^{\g}$, i.e., $\s_X(0,\d)$ and $\d^{\g}$ are of the same order as $\d\downarrow 0$;

(C2)\quad there exist a constant $\kap>0$ such that
\[
\L(\d):=\sup_{\varphi(\d)\ls t\ls T-\d} \sup_{0<h\ls \d}\bigg\vert\frac{\s_X(t,t+h)}{\kap h^{\g}}- 1\bigg\vert\longrightarrow 0\qquad\mbox{as}\ \d\downarrow 0
\]
for every function $\varphi\in\widetilde \Psi$.

If for some constant $\g\in(0,1)$ the second order stochastic process $X$ satisfies  conditions $(C1)$ and $(C2)$, then the Orey index is equal to $\g$ (see \cite{kk1}).

Recently much attention  has been given to studies and applications of Gaussian  processes such as  fractional Brownian motion (fBm), sub-fractional Brownian motion (subfBm), bifractional Brownian motion (bifBm), fractional Ornstein-Uhlenbeck process. All of them are Gaussian processes and they have the Orey indexes. Consequently, examining Gaussian processes with the Orey index we thus examine the processes listed above.

Many authors (see \cite{glad}, \cite{gl}, \cite{IL}, \cite{Cohen-98}, \cite{Coeurjolly-01}, \cite{Begyn-06}, \cite{begyn2},  \cite{rn1},  \cite{ma}) considered an almost sure convergence and asymptotic normality of the generalized quadratic variations associated to the filter $a$ (see \cite{IL})  of a wide class of processes with Gaussian increments.
The strong consistency of the Orey index estimator was proven in \cite{kk1}.

In the papers of Breton et al. \cite{bnp} and Breton and Coeurjolly \cite{bc}, an exact (non-asymptotic) confidence interval for the Hurst index of fBm was derived with the aid of concentration inequalities for quadratic forms of Gaussian process. The obtained confidence intervals for the Hurst parameter were based on a single observation of a discretized sample path of the interval $[0, 1]$ of fBm. Exact confidence intervals for sub-fractional Brownian motion were considered in \cite{liu} but are not sufficiently precise.

The purpose of article is to extend the results of Breton et al. \cite{bnp} and Breton and Coeurjolly \cite{bc} as well as to apply them to Gaussian processes  with the Orey index which may not have stationary increments.

The paper is organized in the following way. In Section 2 we give exact confidence intervals for the Orey index of Gaussian process. Section 3 contains some application results for known Gaussian processes which may not have stationary increments. Finally, in Section 4 some simulations are given in order to illustrate the obtained results. In addition, Appendix includes the R code listings of simulations.

\section{Confidence intervals}

First, we formulate a concentration inequality for a family of Gaussian r.v.'s.
Consider a finite centered Gaussian family $X = \{X_k : k=1,...,M\}$, and write $d_{kj} = \E(X_k X_j)$. Define two quadratic forms associated with $X$ and with some real coefficient $c$:
\begin{equation}
Q_1(c,X) = c\sum_{k=1}^M (X^2_k - d_{kk}),\qquad Q_2(c,X) = 2c^2\sum_{k,j=1}^M X_k X_j d_{kj}.
\end{equation}
The following statement characterizes the tail behavior of $Q_1(c,X)$.
\begin{thm}{$ (\cite{bnp},\cite{bc}) $} Suppose that $Q_1(c,X)$ is not a.s. zero and fix $\alpha\geq 0$ and $\beta>0$.
Assume that $Q_2(c,X)\leq \alpha Q_1(c,X)+\beta$, a.s.-$P$. Then, for all $z>0$, we have
\begin{align*}
\pr(Q_1(c,X)\geq z)\ls& \varphi_r(z;\a,\b):=e^{-z/\a}\Big(1+\frac{\a z}{\b}\Big)^{{\b}/{\a^2}}\\
\pr(Q_1(c,X)\leq -z)\ls& \varphi_l(z;\a,\b):=e^{z/\a}\Big(1-\frac{\a z}{\b}\Big)^{{\b}/{\a^2}}{\bf 1}_{[0,\b/\a]}(z).
\end{align*}
\end{thm}

\begin{rem}{\cite{bc}}\label{bijection} Note that $\varphi_r(\cdot;\a,\b)$ (resp. $\varphi_l(\cdot;\a,\b)$) is a bijective function from $(0,+\infty)$ (resp. $(0,\b/\a)$) to $(0,1)$.
\end{rem}

Next, we apply the obtained concentration inequality to second order quadratic variations.

Let $X$ be a a centered Gaussian process satisfying conditions $(C1)$ and $(C2)$ with the Orey index $\g\in(0,1)$. Denote $Y=(Y_{1,n},\ldots,Y_{n-1,n})$, where
\[
Y_{k,n}=\frac{n^{\g}}{T^\g\varkappa \sqrt{4-2^{2\g}}}\,\D^{(2)}_{k,n}X,\qquad \D^{(2)}_{k,n} X=X\big(t^n_{k+1}\big)-2X\big(t^n_{k}\big) +X\big(t^n_{k-1}\big),
\]
$t^n_k=\frac{kT}{n}$ and $\kap$ is a constant defined in condition $(C2)$. Set $d^{Y,n}_{jk}=\E Y_{j,n}Y_{k,n}$.

\begin{prop}\label{main} Assume that there exists a sequence of real numbers $(\eps_n)$ not depending on $\g$ and such that
\begin{equation}\label{kvadrvar}
\bigg\vert\frac{1}{n-1}\sum_{k=1}^{n-1} d^{Y,n}_{kk}-1\bigg\vert\ls \eps_n, \qquad \eps_n\downarrow 0.
\end{equation}
Then for all $z>0$ we have
\begin{align*}
\pr\bigg(\frac{1}{\sqrt{n-1}}\sum_{k=1}^{n-1}\big[\big(Y_{k,n}\big)^2-d^{Y,n}_{kk}\big]\gs z\bigg)\ls& \varphi_{r,n}(z;\nu_n,\eps_n),\\
\pr\bigg(\frac{1}{\sqrt{n-1}}\sum_{k=1}^{n-1}\big[\big(Y_{k,n}\big)^2-d^{Y,n}_{kk}\big]\ls -z\bigg)\ls& \varphi_{l,n}(z;\nu_n,\eps_n),
\end{align*}
where
\begin{align*}
\varphi_{r,n}(z;\nu_n,\eps_n):=&\exp\bigg\{-\frac{z\sqrt{n-1}}{2\nu_n}\bigg\}\bigg(1+\frac{z}{(\eps_n+1)\sqrt{n-1}}\bigg)^{\frac{(\eps_n+1)(n-1)}{2\nu_n}}\\
\varphi_{l,n}(z;\nu_n,\eps_n)
:=&\exp\bigg\{\frac{z\sqrt{n-1}}{2\nu_n}\bigg\}\bigg(1-\frac{z}{(\eps_n+1)\sqrt{n-1}}\bigg)^{\frac{(\eps_n+1)(n-1)}{2\nu_n}}{\bf 1}_{[0,(\eps_n+1)\sqrt{n-1}\,]}(z),\\
\nu_n=&\max_{1\ls k\ls n-1}\sum_{j=1}^{n-1}  \vert d^{Y,n}_{kj}\vert.
\end{align*}
\end{prop}
\proof Denote
\[
Q_1((n-1)^{-1/2},Y)=\frac{1}{\sqrt{n-1}}\sum_{k=1}^{n-1}\big[\big(Y_{k,n}\big)^2-d^{Y,n}_{kk}\big]
\]
and
\[
Q_2((n-1)^{-1/2},Y)=\frac{2}{n-1}\sum_{k,j=1}^{n-1} Y_k Y_j\, d^{Y,n}_{kj}.
\]
Then
\begin{align*}
&Q_2((n-1)^{-1/2},Y)\\
&\quad
\ls \frac{2}{n-1}\sum_{k,j=1}^{n-1} \vert Y_k\vert\cdot \vert Y_j\vert\cdot \vert d^{Y,n}_{kj}\vert\ls \frac{1}{n-1}\sum_{k,j=1}^{n-1} [Y_k^2+Y_j^2] \cdot \vert d^{Y,n}_{kj}\vert\\
&\quad=\frac{2}{n-1}\sum_{k,j=1}^{n-1} Y_k^2  \vert d^{Y,n}_{kj}\vert=\frac{2}{n-1}\sum_{k=1}^{n-1}Y_k^2\sum_{j=1}^{n-1}   \vert d^{Y,n}_{kj}\vert\\
&\quad\ls \frac{2}{n-1}\sum_{k=1}^{n-1}Y_k^2 \bigg(\max_{1\ls k\ls n-1}\sum_{j=1}^{n-1}  \vert d^{Y,n}_{kj}\vert\bigg)\\
&\quad= \frac{2\nu_n}{\sqrt{n-1}}\bigg(Q_1((n-1)^{-1/2},Y) +\frac{1}{\sqrt{n-1}}\sum_{k=1}^{n-1}\big[d^{Y,n}_{kk}-1\big]+\sqrt{n-1}\bigg)\\
&\quad\ls\alpha_n Q_1((n-1)^{-1/2},Y)+\beta_n,
\end{align*}
where
\[
\a_n=\frac{2\nu_n}{\sqrt{n-1}}\,,\quad \b_n=  2\nu_n (\eps_n+1).
\]
Thus
\begin{align}
&\pr\big(Q_1((n-1)^{-1/2},Y)\gs z\big)\ls \varphi_{r,n}(z;\nu_n,\eps_n)\label{ineq1}\\
&\pr\big(Q_1((n-1)^{-1/2},Y)\ls -z\big)\ls \varphi_{l,n}(z;\nu_n,\eps_n)\label{ineq2}
\end{align}
and the proof is completed.\endproof

\begin{rem} Since $\varphi_{l,n}(z;\cdot,\eps_n)$ and $\varphi_{r,n}(z;\cdot,\eps_n)$ are non-decreasing, inequalities (\ref{ineq1}) and (\ref{ineq2}) remain true with $\nu_n$ replacing $\upsilon_n$, where $\nu_n\ls \upsilon_n$.
\end{rem}

For any $\a\in(0,1)$ and $\upsilon_n$ ($\nu_n\ls \upsilon_n$), denote by $q_{l,n}(\a):=(\varphi_{l,n})^{-1}(\a;\upsilon_n,\eps_n)$ and $q_{r,n}(\a):=(\varphi_{r,n})^{-1}(\a;\upsilon_n,\eps_n)$. For convenience we define
\[
x_{l,n-1}(\a):=1-\frac{q_{l,n-1}(\a)}{\sqrt{n-1}}\quad\mbox{and}\quad x_{r,n-1}(\a):=1+\frac{q_{r,n-1}(\a)}{\sqrt{n-1}}\,.
\]
Note that Remark \ref{bijection} above ensures that for any $\a\in(0,1)$ and for all $n>1$, we have $x_{l,n-1}(\a)>0$.
Set
\[
S_n:=\frac{1}{n-1}\sum_{k=1}^{n-1}\big(\D^{(2)}_{k,n} X \big)^2,\qquad g_{n,T}(x):=2x\ln(n/T)-\ln (4-2^{2x}),\qquad x\in (0,1),
\]
and
\[
\ln^*x=\begin{cases} -\infty &\mbox{if}\quad x= 0,\\
\ln x&\mbox{if}\quad x>0.
\end{cases}
\]
The function $g_{n,T}(x)$ is a strictly increasing bijection from $(0,1)$ to $(-\ln 3,+\infty)$ if $n>T$.

\begin{thm}\label{pagrindine} Let $\a\in(0,1)$. Assume that conditions of Proposition \ref{main} are satisfied and there are constants $\upsilon_n$ such that $\nu_n\ls \upsilon_n$. Then
\[
\mathbb{P}\big(\g\in\big[\g_n^{inf}(\a),\g_n^{sup}(\a)]\big)\gs 1-\a,
\]
where
\begin{align*}
\g_n^{inf}(\a):=&\max\bigg(0,g^{-1}_{n,T}\bigg[\max\bigg(\ln^*\bigg(\bigg(\frac{x_{l,n-1}(\a/2) -\eps_n}{S_n}\,\kap^2\bigg)\lor 0\bigg),-\ln 3\bigg)\bigg],\\
\g_n^{sup}(\a):=&\min\bigg(1,g^{-1}_{n,T}\bigg[\ln\bigg(\frac{x_{r,n-1}(\a/2)+\eps_n}{S_n}\,\kap^2\bigg)\bigg]\bigg).
\end{align*}
\end{thm}
\proof Denote
\[
Z_n=(n-1)^{-1/2}V_n^{(2)}(Y,2)-\sqrt{n-1}=\sqrt{n-1} \big[(n-1)^{-1}V_n^{(2)}(Y,2)-1\big],
\]
where
\[
V_n^{(2)}(Y,2)=\sum_{i=1}^{n-1}\big(\D^{(2)}_{k,n} Y \big)^2.
\]
Then
\begin{align*}
&\bigg\{-q_{l,n-1}(\a/2)\ls(n-1)^{-1/2}\sum_{k=1}^{n-1}\big[\big(Y_{k,n}\big)^2-d^{Y,n}_{kk}\big]\ls q_{r,n-1}(\a/2)\bigg\}\\
&\quad=\bigg\{-q_{l,n-1}(\a/2)+(n-1)^{-1/2}\sum_{k=1}^{n-1}\big[d^{Y,n}_{kk}-1\big]\ls Z_n \ls q_{r,n-1}(\a/2)+(n-1)^{-1/2}\sum_{k=1}^{n-1}\big[d^{Y,n}_{kk}-1\big]\bigg\}\\
&\quad=\bigg\{x_{l,n-1}(\a/2)+(n-1)^{-1}\sum_{k=1}^{n-1}\big[d^{Y,n}_{kk}-1\big]\ls \frac{n^{2\g}}{T^{2\g}\kap^2(4-4^\g)}\,S_n\\
&\qquad\qquad\ls x_{r,n-1}(\a/2)+(n-1)^{-1}\sum_{k=1}^{n-1}\big[d^{Y,n}_{kk}-1\big]\bigg\}\\
&\quad\subset \bigg\{x_{l,n-1}(\a/2)-\eps_n\ls \frac{n^{2\g}}{T^{2\g}\kap^2(4-2^{2\g})}\,S_n\ls x_{r,n-1}(\a/2)+\eps_n\bigg\}\\
&\quad=\bigg\{\big(x_{l,n-1}(\a/2)-\eps_n\big)\lor 0\ls \frac{n^{2\g}}{T^{2\g}\kap^2(4-2^{2\g})}\,S_n\ls x_{r,n-1}(\a/2)+\eps_n\bigg\}\\
&\quad=\big\{\ln^*\big(\big(x_{l,n-1}(\a/2)-\eps_n\big)\lor 0\big)-\ln S_n+\ln \kap^2\ls g_{n,T}(\g)\ls \ln(x_{r,n-1}(\a/2)+\eps_n)-\ln S_n+\ln \kap^2\big\}\\
&\quad =\bigg\{\bigg[\ln^*\bigg(\frac{(x_{l,n-1}(\a/2)-\eps_n)\lor 0}{S_n}\,\kap^2\bigg)\bigg]\lor(-\ln 3)\ls g_{n,T}(\g)\ls \ln\bigg(\frac{x_{r,n-1}(\a/2)+\eps_n}{S_n}\,\kap^2\bigg)\bigg\}.
\end{align*}
Note that
\[
\pr\bigg(-q_{l,n-1}(\a/2)\ls \frac{1}{\sqrt{n-1}} \sum_{k=1}^{n-1}\big[\big(Y_{k,n}\big)^2-d^{Y,n}_{kk}\big]\ls q_{r,n-1}(\a/2)\bigg)\gs 1-\a.
\]
Thus
\[
\pr\big(\g\in\big[\g_n^{inf}(\a),\g_n^{sup}(\a)]\big)\gs 1-\a.
\]
The proof is completed.

\section{Applications}

In this section we obtain the confidence intervals for subfBm, bifBm and the fractional Ornstein-Uhlenbeck process. For this purpose we apply the Theorem \ref{pagrindine}. In order to apply the Theorem \ref{pagrindine}, it suffices to find the sequence of real numbers $(\eps_n)$ in the estimation (\ref{kvadrvar}) and estimate $\nu_n$. In the considered cases as the special case appears the Brownian motion. We exclude it from consideration in view of its properties (in particular, independent increaments). It is easy to see that
\[
\frac{n}{2T}\, d^{B,n}_{kk}=1,\qquad \frac{n}{2T}\,\max_{1\ls k\ls n-1}\sum_{j=1}^{n-1}  \vert d^{B,n}_{jk}\vert\ls 2,
\]
where $d^{B,n}_{jk}=\E \D^{(2)}_{j,n} B \D^{(2)}_{k,n} B$ and $B$ is Brownian motion.

\subsection{Sub-fractional Brownian motion}

\begin{dfn}{\rm (\cite{BGT})}
A \textbf{sub-fractional Brownian motion} (subfBm) with the index $H$, $H\in(0,1)$, is a mean zero Gaussian stochastic process $S^H=(S^H_t, t\gs 0)$ with the covariance function
\[
G_H(s,t):=s^{2H} +t^{2H}-\frac{1}{2}\big[(s+t)^{2H}+\vert s-t\vert^{2H}\big].
\]
\end{dfn}
The case $H =1/2$ corresponds to the Brownian motion. For $H\neq 1/2$ this process has some of the main properties of fBm, but its increments are not stationary.

The incremental variance function of subfBm is of the following form
\begin{equation}\label{increm1}
\s_{S^H}^2(s,t)=\E\vert S^H_t-S^H_s\vert^2= \vert t-s\vert^{2H}+(s+t)^{2H}-2^{2H-1}(t^{2H}+s^{2H}).
\end{equation}
For any $0\ls s\ls t\ls T$ the inequalities (see \cite{BGT})
\begin{align}
&(t-s)^{2H}\ls\s_{S^H}^2(s,t)\ls (2-2^{2H-1})(t-s)^{2H}, \qquad\mbox{if}\quad 0<H<1/2,\label{subf1}\\
&(2-2^{2H-1})(t-s)^{2H}\ls\s_{S^H}^2(s,t)\ls(t-s)^{2H}, \qquad\mbox{if}\quad 1/2<H<1\label{subf2}
\end{align}
hold.

It is known  that for subfBm the Orey index is equal to $H$ (see \cite{kk1}). Now we prove the following lemma.

\begin{lem}\label{sublema} Assume that $S^H=\{S^H(t): t\in[0,T]\}$ is a subfBm. If $H\neq 1/2$ then \[
\bigg\vert\frac{1}{n-1}\sum_{k=1}^{n-1} d^{S^H,n}_{kk}-1\bigg\vert\ls \bigg(\frac{T}{n}\bigg)^{2/3}\bigg(\frac{1}{6T}\,\frac{n}{n-1}+\frac{33}{9\ln 4}\bigg)\,,
\]
where
\[
d^{S^H,n}_{kk}=\frac{n^{2H}}{T^{2H} (4-2^{2H})}\,\E \big(\D^{(2)}_{k,n} S^H\big)^2.
\]
\end{lem}
\proof Observe that the following equality
\begin{align*}
{\bf E}\big(S^H_{t+h}-2S^H_t+S^H_{t-h}\big)^2
=& (4-2^{2H})h^{2H}-2^{2H-1}(t+h)^{2H}-3\cdot 2^{2H}t^{2H} \\
&-2^{2H-1}(t-h)^{2H}+2(2t+h)^{2H}+2(2t-h)^{2H}
\end{align*}
holds. Thus
\begin{align*}
d^{S^H,n}_{kk}=&1-\frac{2^{2H-1}(k+1)^{2H}+3\cdot 2^{2H}k^{2H}+2^{2H-1}(k-1)^{2H}-2(2k+1)^{2H}-2(2k-1)^{2H}}{4-2^{2H}}\\
=&1-\frac{b(k,H)}{4-2^{2H}}\,.
\end{align*}
For simplicity we shall omit the index $n$ for $d^{S^H,n}_{kk}$. Using computer modeling we obtain the inequalities
\begin{equation}\label{kovariance}
\max_{H\in(0,1)}  d^{S^H}_{kk}\ls \max_{H\in(0,1)}  d^{S^H}_{11}=\frac 76\quad \mbox{and}\quad \max_{H\in(0,1)}\frac{\vert b(k,H)\vert}{4-2^{2H}} \ls \max_{H\in(0,1)}\frac{\vert b(1,H)\vert}{4-2^{2H}}\ls \frac{1}{6}\,.
\end{equation}
Let $\varphi\in\widetilde \Psi$ and denote $\tau_n:=[\varphi(Tn^{-1})nT^{-1}]$, where $[a]$ is an integer part of $a$. Then
\begin{align*}
\bigg\vert\frac{1}{n-1}\sum_{k=1}^{n-1} d^{S^H}_{kk}-1\bigg\vert
\ls& \frac{1}{n-1}\sum_{k=1}^{\tau_n}\frac{\vert b(k,H)\vert}{4-2^{2H}}+\max_{\tau_n+1\ls k\ls n-1}\frac{\vert b(k,H)\vert}{4-2^{2H}}\\
=&I_1+I_2.
\end{align*}
It is clear that
\[
I_1\ls \frac{\tau_n}{n-1}\,\max_{1\ls k\ls \tau_n} \frac{\vert b(k,H)\vert}{4-2^{2H}}\ls \frac{\tau_n}{n-1}\,\frac{\vert b(1,H)\vert}{4-2^{2H}}\ls \frac{1}{6}\,\varphi(Tn^{-1})T^{-1}\,\frac{n}{n-1}\,.
\]
Now we estimate $b(k,H)$ for $k\gs 2$. Using the formula
\[
(1+x)^\a=1+\sum_{k=1}^\infty\frac{\a(\a-1)\cdots(\a-k+1)}{k!}\,x^k\qquad\mbox{for}\ -1<x<1,
\]
we obtain
\begin{align*}
b(k,H)=&2^{2H-1}(k+1)^{2H}+3\cdot 2^{2H}k^{2H} +2^{2H-1}(k-1)^{2H}-2(2k+1)^{2H}-2(2k-1)^{2H}\\
=&2^{2H-1}k^{2H}\bigg[\bigg(1+\frac{1}{k}\bigg)^{2H}+6  +\bigg(1-\frac{1}{k}\bigg)^{2H}-4\bigg(1+\frac{1}{2k}\bigg)^{2H}-4\bigg(1-\frac{1}{2k}\bigg)^{2H}\bigg]\\
=&2^{2H-1}k^{2H}\bigg[\sum_{m=1}^\infty \frac{2H(2H-1)\cdots(2H-m+1)}{m!}\,\Big(\frac 1k\Big)^m \big[1+(-1)^m-2^{2-m}-(-1)^m2^{2-m}\big]\bigg]\\
=&2^{2H}k^{2H}\bigg[\sum_{m=2}^\infty \frac{2H(2H-1)\cdots(2H-2m+1)}{(2m)!}\,\Big(\frac 1k\Big)^{2m} \big(1-2^{2-2m}\big)\bigg].
\end{align*}
Note that the sign of $2H(2H-1)(2H-2)\cdots(2H-2m+1)$ is the same as that of $2H-1$. Thus
\begin{align*}
\frac{\vert b(k,H)\vert}{4-2^{2H}}\ls&\frac{2^{2H}k^{2H}}{4-2^{2H}}\bigg[\sum_{m=2}^\infty \frac{ 2H\vert 2H-1\vert(2-2H)\cdots(2m-1-2H)}{(2m)!}\,\Big(\frac 1k\Big)^{2m} \big[1-2^{2-2m}\big]\bigg]\\
\ls&\frac{3\cdot 2^{2H}k^{2H}}{4\ln 4}\sum_{m=2}^\infty \frac{ 2(2m-1)!}{(2m)!}\,\Big(\frac{1}{k}\Big)^{2m} \\
\ls& \frac{3\cdot 2^{2H}}{4\ln 4}\,k^{2H}\sum_{m=2}^\infty \frac{1}{m}\,\Big(\frac{1}{k^2}\Big)^m   \ls \frac{3\cdot 2^{2H}}{4\ln 4}\,k^{2H-2}\sum_{m=2}^\infty \frac{1}{m-1}\,\Big(\frac{1}{k^2}\Big)^{m-1}\\
=& -\frac{3\cdot 2^{2H}}{4\ln 4}\,k^{2H-2}\ln\Big(1-\frac {1}{k^2}\Big)\ls \frac{33\cdot 2^{2H}}{36\ln 4}\,k^{2H-4}\ls \frac{33}{9\ln 4}\,k^{-2}
\end{align*}
since
\begin{align*}
&2H\vert 1-2H\vert (2-2H)\cdots (2m-1-2H)\ls 2(1-H)(2m-1)!\,,\\ &\log(1-u)=-\sum_{k=1}^{\infty} \frac{u^k}{k}\quad\mbox{if}\ 0\leq u<1\qquad\mbox{and}\qquad -\log(1-u)\ls \frac{11}{9}u \quad\mbox{if}\ 0\leq u \ls 1/4,\\
&(1-H)\ln 4<4-2^{2H}.
\end{align*}
So
\begin{align*}
\bigg\vert\frac{1}{n-1}\sum_{k=1}^{n-1} d^{S^H}_{kk}-1\bigg\vert\ls& \frac{1}{6}\,\varphi(Tn^{-1})T^{-1}\,\frac{n}{n-1}+\frac{33}{9\ln 4}\, \frac{1}{(\tau_n+1)^2}\\
\ls&\frac{1}{6}\,\varphi(Tn^{-1})T^{-1}\,\frac{n}{n-1}+\frac{33}{9\ln 4}\, \frac{T^2}{\varphi^2(Tn^{-1})n^2}\,.
\end{align*}
Let $\varphi(Tn^{-1})=(Tn^{-1})^{2/3}$. It belongs to the class of functions $\widetilde \Psi$. After putting into the obtained inequality we get
\[
\bigg\vert\frac{1}{n-1}\sum_{k=1}^{n-1} d^{S^H}_{kk}-1\bigg\vert
\ls \frac{1}{6}\, \frac{1}{T^{1/3} n^{2/3}}\,\frac{n}{n-1}+\frac{33}{9\ln 4}\bigg(\frac{T}{n}\bigg)^{2/3}
=\bigg(\frac{T}{n}\bigg)^{2/3}\bigg(\frac{2}{15 T}\,\frac{n}{n-1}+\frac{33}{9\ln 4}\bigg).
\]

\begin{lem}
Assume that $S^H=\{S^H_t: t\in[0,T]\}$ is a subfBm. If $H\neq 1/2$ then
\begin{equation}\label{sfBm2}
\max_{1\ls k\ls n-1}\sum_{j=1}^{n-1}  \vert d^{S^H,n}_{jk}\vert\ls \frac{9}{2}\,,
\end{equation}
where
\[
d^{S^H,n}_{jk}=\frac{n^{2H}}{T^{2H} (4-2^{2H})}\,\widehat d_{jk}^{S^H,n}, \qquad\widehat d_{jk}^{S^H,n}=\E \D^{(2)}_{j,n} S^H \D^{(2)}_{k,n} S^H.
\]
\end{lem}
\proof The fourth order mixed partial derivative of the covariance function $G_H(s,t)$ is of the following form
\[
\frac{\partial^4 G_H}{\partial s^2\partial t^2}(s,t)
= -C_H \bigg[\frac{1}{\vert s-t\vert^{2(2-H)}}+\frac{1}{(s+t)^{2(2-H)}} \bigg]
\]
for each $s,t>0$ such that $s\neq t$, where $C_H=H(2H-1)(2H-2)(2H-3)$. Since the covariance function $G_H(s,t)$ is continuous in $[0,T]^2$ and the derivative $\frac{\partial^4R}{\partial s^2\partial t^2}$ is continuous in $(0,T]^2/\{s=t\}$ then for $H\neq 1/2 $ and $j\neq 1$ or $k\neq 1$
\begin{align*}
\widehat d_{jk}^{S^H,n}=& \int_{t^n_j}^{t^n_{j+1}}du \int^u_{u-T/n}dv \int_{t^n_k}^{t^n_{k+1}}dx\int^x_{x-T/n}  \frac{\partial^4 G_H}{\partial s^2\partial t^2}(v,y)\,dy.
\end{align*}
Assume that $R_H$ is the covariance function of the fBm $B^H$. Then the derivative
\[
\frac{\partial^4 R_H}{\partial s^2\partial t^2}(s,t)
= -\frac{C_H}{\vert s-t\vert^{4-2H}}
\]
of the covariance function $R_H$ is continuous in $(0,T]^2/\{s=t\}$ and
\begin{align*}
\widehat d_{jk}^{B^H,n}:=&\E (\D^{(2)}_{k,n} B^H)^2=-\int_{t^n_j}^{t^n_{j+1}}du \int^u_{u-T/n}dv \int_{t^n_k}^{t^n_{k+1}}dx\int^x_{x-T/n}  \frac{C_H }{(v-y)^{4-2H}}\,dy \\ =&\rho_H(j-k)\,\frac{T^{2H}}{n^{2H}}\,,
\end{align*}
where
\[
\rho_H(r)=\frac12\big(-|r-2|^{2H}+4|r-1|^{2H}-6|r|^{2H}+4|r+1|^{2H}-|r+2|^{2H}\big), \quad r\in \mathbb{N}\cup\{0\}.
\]
Note that $(s+t)\gs \vert s-t\vert$  and
\begin{align*}
\big\vert \widehat d_{jk}^{S^H,n}-\widehat d_{jk}^{B^H,n}\big\vert=&\bigg\vert\int_{t^n_j}^{t^n_{j+1}}du \int^u_{u-T/n}dv \int_{t^n_k}^{t^n_{k+1}}dx\int^x_{x-T/n}  \bigg(\frac{\partial^4 G_H}{\partial s^2\partial t^2}(v,y)+\frac{C_H}{(v-y)^{4-2H}}\bigg)\,dy\bigg\vert\\
=&\vert C_H\vert\,\bigg\vert\int_{t^n_j}^{t^n_{j+1}}du \int^u_{u-T/n}dv \int_{t^n_k}^{t^n_{k+1}}dx\int^x_{x-T/n}  \frac{dy}{(v+y)^{4-2H}}\bigg\vert\\
\ls&\vert C_H\vert\,\bigg\vert\int_{t^n_j}^{t^n_{j+1}}du \int^u_{u-T/n}dv \int_{t^n_k}^{t^n_{k+1}}dx\int^x_{x-T/n}  \frac{dy}{(v-y)^{4-2H}}\bigg\vert\\
=& \vert\rho_H(j-k)\vert(Tn^{-1})^{2H}
\end{align*}
for $\vert j-k\vert\gs 1$. Thus
\[
\big\vert \widehat d_{j,k}^{S^H,n}\big\vert \ls 2 \vert\rho_H(j-k)\vert(Tn^{-1})^{2H}\qquad\mbox{for}\ \vert j-k\vert\gs 1.
\]
It still remains to prove the cases when $j=1$ and $k>1$ or $k=1$ and $j>1$.
Set $t^{n,\eps}_2=t^n_2-T\eps/n$ and  $u^\eps=u+T\eps/n$. Since
\begin{align*}
\big\vert \widehat d_{1k}^{S^H,n}-\widehat d_{1k}^{B^H,n}\big\vert=&\lim_{\eps\to 0}\bigg\vert\int_{t^n_1}^{t^{n,\eps}_2}du \int^u_{u^\eps-T/n}dv \int_{t^n_k}^{t^n_{k+1}}dx\int^x_{x-T/n}  \bigg(\frac{\partial^4 G_H}{\partial s^2\partial t^2}(v,y)+\frac{C_H}{(v-y)^{4-2H}}\bigg)\,dy\bigg\vert\\
\ls& \vert C_H\vert\,\lim_{\eps\to 0}\bigg\vert\int_{t^n_1}^{t^{n,\eps}_2}du \int^u_{u^\eps-T/n}dv \int_{t^n_k}^{t^n_{k+1}}dx\int^x_{x-T/n}  \frac{dy}{(v-y)^{4-2H}}\bigg\vert\\
=&\frac 12(Tn^{-1})^{2H}\, \lim_{\eps\to 0}\big\vert[(k+1)-(2-\eps)]^{2H}-2\cdot k^{2H}+(k+1-\eps)^{2H}-2[k-(2-\eps)]^{2H}\\
&+4(k-1)^{2H}-2(k-\eps)^{2H}+[(k-1)-(2-\eps)]^{2H}-2(k-2)^{2H}+(k-1-\eps)^{2H}\big\vert\\
=&(Tn^{-1})^{2H}\,\vert\rho_H(k-1)\vert
\end{align*}
then the inequality
\[
\big\vert \widehat d_{1,k}^{S^H}\big\vert\ls  2 \vert\rho_H(k-1)\vert(Tn^{-1})^{2H}
\]
holds.  A similar argument yields
\[
\big\vert \widehat d_{j,1}^{S^H}\big\vert\ls  2 \vert\rho_H(j-1)\vert(Tn^{-1})^{2H}.
\]

Now we shall prove the statement of the lemma. We will use the estimate (\ref{kovariance}) and the equalities $\rho_H(-r)=\rho_H(r)$, $\rho_H(0)=4-2^{2H}$, and $\rho_H(1)=-\frac 12(7-4\cdot 2^{2H}+3^{2H})$.  Note that for $H\neq 1/2$
\begin{align*}
&\max_{1\ls k\ls n-1}\sum_{j=1}^{n-1}  \big\vert d^{S^H}_{jk}\big\vert\ls\frac 76+2\max_{1\ls k\ls n-1}\sum_{j=1\atop j\neq k}^{n-1}\frac{\vert \rho_H(j-k)\vert}{\vert \rho_H(0)\vert}\\
&\quad\ls\frac 76+4\sum_{j=1}^\infty\frac{\vert \rho_H(j)\vert}{\vert \rho_H(0)\vert}=\frac 76+2\,\frac{7-4\cdot 2^{2H}+3^{2H}}{4-4^H}+4 \sum_{j=2}^\infty\frac{\vert \rho_H(j)\vert}{\vert \rho_H(0)\vert}\,.
\end{align*}
In \cite{bc} it was proven that
\[
\sum_{j=2}^\infty\frac{\vert \rho_H(j)\vert}{\vert \rho_H(0)\vert}=\frac 12\cdot\sgn(2H-1)\,\frac{3-3\cdot 2^{2H}+3^{2H}}{4-4^H}\,.
\]
Thus
\begin{align*}
\max_{1\ls k\ls n-1}\sum_{j=1}^{n-1}  \big\vert d^{S^H}_{jk}\big\vert
\ls&\frac 76+2\,\frac{7-4\cdot 2^{2H}+3^{2H}}{4-4^H}-2\cdot\sgn(2H-1)\,\frac{3-3\cdot 2^{2H}+3^{2H}}{4-4^H}\\
=&\begin{cases} \frac 76+2\,\frac{10-7\cdot 4^H+2\cdot 3^{2H}}{4-4^H} &\mbox{for}\quad H<1/2,\\
\frac 76+2&\mbox{for}\quad 1/2<H<1
\end{cases}\\
\ls&\begin{cases} \frac 76+\frac{10}{3} &\mbox{for}\quad H<1/2,\\
\frac 76+2&\mbox{for}\quad 1/2<H<1.
\end{cases}
\end{align*}

\subsection{Bifractional Brownian motion}

\begin{dfn}{\rm (\cite{HV})}
A \textbf{bifractional Brownian motion}  (bifBm)  $B^{KH}=(B^{KH}_t, t\gs 0)$ with parameters $H\in(0,1)$ and $K\in(0,1]$ is a centered Gaussian process with the covariance function
\[
F_{KH}(t,s)=2^{-K}\big((t^{2H} +s^{2H})^K-\vert t-s\vert^{2HK}\big),\qquad s,t\gs 0.
\]
\end{dfn}
The incremental variance function of bifBm is
\[
\s_{B^{H,K}}^2(s,t)=\E\vert B^{H,K}_t-B^{H,K}_s\vert^2= 2^{1-K}\big[\vert t-s\vert^{2HK}-(t^{2H}+s^{2H})^K\big]+t^{2HK}+s^{2HK}.
\]
Let $H\in(0,1)$ and $K\in(0,1]$. Then
\begin{equation}\label{nelyg9}
2^{-K}\vert t-s\vert^{2HK}\ls\s_{B^{H,K}}^2(s,t)\ls 2^{1-K}\vert t-s\vert^{2HK}
\end{equation}
for all $s,t\in[0,\infty)$ (see \cite{HV}).

It is known  that for bifBm the Orey index is equal to $HK$ (see \cite{kk1}). If $K=1$ then bifBm becomes fBm, hence we ignore this case.

\begin{lem} Assume that $B^{H,K}=\{B^{H,K}(t): t\in[0,T]\}$ is a bifBm with $K\in(0,1)$ and $H\in(0,1/2)$. Then
\[
\bigg\vert\frac{1}{n-1}\sum_{k=1}^{n-1} d^{B^{KH}}_{kk}-1\bigg\vert\ls \bigg(\frac{T}{n}\bigg)^{1/2}\bigg(\frac{1}{6T}\,\frac{n}{n-1}+\frac{22}{ 9\ln 4}\bigg),
\]
where
\[
d^{B^{KH}}_{kk}=\frac{n^{2H}}{2^{1-K}T^{2H} (4-2^{2H})}\,\E \big(\D^{(2)}_{k,n} B^{H,K}\big)^2.
\]
\end{lem}
\begin{rem} Without this restriction for $H$ the expressions becomes more complicated.
\end{rem}
\proof The proof of the lemma follows the outlines of the proof of Lemma \ref{sublema}. Observe that the following equality
\[
d^{B^{KH}}_{kk}=1-\frac{b(k,H,K)}{4-4^{KH}}
\]
holds, where
\begin{align*}
b(k,H,K)=&2[(k+1)^{2H}+k^{2H}]^K+2[k^{2H}+(k-1)^{2H}]^K\\
&-2^{K-1}[(k+1)^{2KH}+4k^{2KH}+(k-1)^{2KH}]-[(k+1)^{2H}+(k-1)^{2H}]^K.
\end{align*}
By computer modeling we obtain inequalities
\begin{equation}\label{covariance1}
\max_{H,K\in(0,1)} d^{B^{KH}}_{kk}\ls 1\quad\mbox{and}\quad \max_{H,K\in(0,1)}\frac{ b(k,H,K)}{4-2^{2KH}}\ls \frac 16
\end{equation}
for all $k\gs 1$. Let $\varphi\in\widetilde \Psi$ and $\tau_n=[\varphi(Tn^{-1})nT^{-1}]$. Then
\begin{align*}
\bigg\vert\frac{1}{n-1}\sum_{k=1}^{n-1} d^{B^{KH}}_{kk}-1\bigg\vert
\ls& \frac{1}{n-1}\sum_{k=1}^{\tau_n}\frac{ b(k,H,K)}{4-2^{2KH}}+\max_{\tau_n+1\ls k\ls n-1}\frac{b(k,H,K)}{4-2^{2KH}}\\
=&I_1+I_2.
\end{align*}
and
\[
I_1\ls \frac{\tau_n}{n-1}\,\max_{1\ls k\ls \tau_n} \frac{ b(k,H,K)}{4-2^{2KH}}\ls \frac 16\,\varphi(Tn^{-1})T^{-1}\,\frac{n}{n-1}\,.
\]
Using the inequality $(a+b)^K\gs 2^{K-1}(a^K+b^K)$, $a,b\gs 0$, $0<K<1$, and the formula
\[
(1+x)^\a=1+\sum_{k=1}^\infty\frac{\a(\a-1)\cdots(\a-k+1)}{k!}\,x^k,\qquad\mbox{for}\ -1<x<1,
\]
we obtain
\begin{align*}
\vert b(k,H,K)\vert
=&k^{2HK}\Bigg\vert 2\bigg[\bigg(1+\frac 1k\bigg)^{2H}+1\bigg]^K-2^{K-1}\bigg[\bigg(1+\frac 1k\bigg)^{2KH}
+4+\bigg(1-\frac 1k\bigg)^{2KH}\bigg]\\
&+2\bigg[1+\bigg(1-\frac 1k\bigg)^{2H}\bigg]^K-\bigg[\bigg(1+\frac 1k\bigg)^{2H}+\bigg(1-\frac 1k\bigg)^{2H}\bigg]^K\Bigg\vert\\
\ls& 2k^{2HK}\Bigg\{2^{1-K}\bigg[\bigg(1+\frac 1k\bigg)^{2H}
+2+\bigg(1-\frac 1k\bigg)^{2H}\bigg]^K\\
&-2^{K-1}\bigg[\bigg(1+\frac 1k\bigg)^{2KH}
+2+\bigg(1-\frac 1k\bigg)^{2KH}\bigg]\Bigg\}\\
=&2k^{2HK}\Bigg\{ 2^{1-K} \bigg[4+\sum_{m=1}^\infty \frac{2H(2H-1)\cdots(2H-2m+1)}{(2m)!}\,\Big(\frac 1k\Big)^{2m}\bigg]^K\\
&-2^{K-1}\bigg[4+\sum_{m=1}^\infty \frac{2KH(2KH-1)\cdots(2KH-2m+1)}{(2m)!}\,\Big(\frac 1k\Big)^{2m}\bigg]\Bigg\}:=\widehat b(k,H,K).
\end{align*}
Note that the sign of $2KH(2KH-1)(2KH-2)\cdots(2KH-2m+1)$ is the same as that of $2KH-1$. Thus the estimate of $b(k,H,K)$ depends on the signs of $2H-1$ and $2KH-1$. Then for $H\in(0,1/2)$ and $k\gs 2$
\begin{align*}
\frac{\widehat b(k,H,K)}{4-4^{KH}}
\ls& \frac{ 2^K k^{2HK}}{4-4^{KH}}\sum_{m=1}^\infty \frac{2KH(1-2KH)\cdots(2m-1-2KH)}{(2m)!}\,\Big(\frac 1k\Big)^{2m}\\
\ls& \frac{ 2^K k^{2HK}(1-KH)}{4-4^{KH}} \sum_{m=1}^\infty \frac{1}{m}\,\Big(\frac 1k\Big)^{2m}\ls \frac{ 2^K k^{2HK}}{\ln 4} \sum_{m=1}^\infty \frac{1}{m}\,\Big(\frac 1k\Big)^{2m}\\
=&-\frac{ 2^K k^{2HK}}{\ln 4}\,\ln\bigg(1- \frac{1}{k^2}\bigg) \ls \frac{11\cdot 2^K}{9\ln 4}\,\frac{1}{k^{2-2KH}}\ls \frac{22}{ 9\ln 4}\,k^{-1}\,.
\end{align*}
Thus
\begin{align*}
\bigg\vert\frac{1}{n-1}\sum_{k=1}^{n-1} d^{B^{KH}}_{kk}-1\bigg\vert
\ls&\frac 16\,\varphi(Tn^{-1})T^{-1}\,\frac{n}{n-1}+\frac{22}{ 9\ln 4}\, \frac{T}{\varphi(Tn^{-1})n}\,.
\end{align*}
Let $\varphi(Tn^{-1})=(Tn^{-1})^{1/2}$. It belongs to the class of functions $\widetilde \Psi$. After plugging it into the obtained inequality we get
\begin{align*}
\bigg\vert\frac{1}{n-1}\sum_{k=1}^{n-1} d^{B^{KH}}_{kk}-1\bigg\vert
\ls&  \frac 16\,\bigg(\frac{T}{n}\bigg)^{1/2}\,\frac{n}{T(n-1)}+\frac{22}{ 9\ln 4}\, \frac{T}{(Tn^{-1})^{1/2} n}\\
\ls& \frac 16\,\bigg(\frac{1}{nT}\bigg)^{1/2}\,\frac{n}{n-1} +\frac{22}{ 9\ln 4}\,\bigg(\frac{T}{n}\bigg)^{1/2}
=\bigg(\frac{T}{n}\bigg)^{1/2}\bigg(\frac{1}{6T}\,\frac{n}{n-1}+\frac{22}{ 9\ln 4}\bigg).
\end{align*}

\begin{lem}
Assume that $B^{H,K}=\{B^{H,K}(t): t\in[0,T]\}$ is a bifBm with $K\in(0,1)$ and $H\in(0,1/2)$. Then
\begin{equation}
\max_{1\ls k\ls n-1}\sum_{j=1}^{n-1}  \vert d^{B^{KH}}_{jk}\vert\ls 5.005.
\end{equation}
\end{lem}
\proof The fourth order mixed partial derivative of the covariance function $R^{KH}(s,t)$ has the form
\begin{align*}
\frac{\partial^4 F_{KH}}{\partial s^2\partial t^2}(s,t)
=&-\frac{\widehat C^{(1)}_{KH}}{\vert s-t\vert^{4-2KH}}
+\widehat C^{(2)}_{KH}\,(st)^{4H-2}\big(s^{2H}+t^{2H}\big)^{K-4}\nonumber\\
&+\widehat C^{(3)}_{KH}(st)^{2H-2}\big(s^{2H}+t^{2H}\big)^{K-2},
\end{align*}
for each $s,\,t>0$ such that $s\neq t$, where
\begin{align*}
\widehat C^{(1)}_{KH}=&2HK(2KH-1)(2HK-2)(2HK-3)2^{-K},\\
\widehat C^{(2)}_{KH}=& K(K-1)(K-2)(K-3)(2H)^4 2^{-K},\\
\widehat C^{(3)}_{KH}=& K(K-1)(2H)^2(2H-1)2^{-K}(2KH-2H-1).
\end{align*}
Since $2s^Ht^H\ls s^{2H}+t^{2H}$ and $K-2<0$, $K-4<0$ it follows that
\begin{align}
(st)^{2H-2}\big(s^{2H}+t^{2H}\big)^{K-2}\ls& 2^{K-2}(st)^{KH-2},\label{ineq3}\\
(st)^{4H-2}\big(s^{2H}+t^{2H}\big)^{K-4}\ls& 2^{K-4}(st)^{KH-2}.\label{ineq4}
\end{align}
Let $\widetilde B^{KH}$ be a fBm with the Orey index $KH$. Assume that $KH\neq 1/2$ and $j\neq 1$ or $k\neq 1$. Using the inequalities (\ref{ineq3}), (\ref{ineq4}) and
\[
vy=[y+(v-y)]y\gs (v-y)y\gs (v-y)\,\frac Tn,\qquad  y\gs t^n_1,
\]
we obtain
\begin{align*}
&\big\vert \widehat d_{jk}^{B^{KH}}-2^{1-K}\widehat d_{jk}^{\widetilde B^{KH}}\big\vert\\
&\quad\ls\vert \widehat C^{(2)}_{KH}\vert \bigg\vert\int_{t^n_j}^{t^n_{j+1}}du \int^u_{u-T/n}dv \int_{t^n_k}^{t^n_{k+1}}dx\int^x_{x-T/n}  \frac{dy}{(vy)^{2-4H}(v^{2H}+y^{2H})^{4-K}} \bigg\vert\\
&\qquad+\vert \widehat C^{(3)}_{KH}\vert \bigg\vert\int_{t^n_j}^{t^n_{j+1}}du \int^u_{u-T/n}dv \int_{t^n_k}^{t^n_{k+1}}dx\int^x_{x-T/n}  \frac{dy}{(vy)^{2-2H}(v^{2H}+y^{2H})^{2-K}} \bigg\vert\\
&\quad\ls\bigg\vert\int_{t^n_j}^{t^n_{j+1}}du \int^u_{u-T/n}dv \int_{t^n_k}^{t^n_{k+1}}dx\int^x_{x-T/n}  \frac{2^{K-2}\vert \widehat C^{(2)}_{KH}\vert+2^{K-4}\vert \widehat C^{(3)}_{KH}\vert}{(vy)^{2-KH}}\,dy\bigg\vert\\
&\quad\ls\bigg(\frac{n}{T}\bigg)^{2-KH}\bigg\vert\int_{t^n_j}^{t^n_{j+1}}du \int^u_{u-T/n}dv \int_{t^n_k}^{t^n_{k+1}}dx\int^x_{x-T/n}  \frac{2^{K-2}\vert \widehat C^{(2)}_{KH}\vert+2^{K-4}\vert \widehat C^{(3)}_{KH}\vert}{(v-y)^{2-KH}}\,dy\bigg\vert\\
&\quad =\bigg(\frac{n}{T}\bigg)^{2-KH}\frac{2^{K-2}\vert \widehat C^{(2)}_{KH}\vert+2^{K-4}\vert \widehat C^{(3)}_{KH}\vert}{\vert C_{KH}\vert}\,\vert \widehat\rho_{KH}(j-k)\vert\,(Tn^{-1})^{KH+2}\\
&\quad=\bigg(\frac{T}{n}\bigg)^{2KH}\frac{2^{K-2}\vert \widehat C^{(2)}_{KH}\vert+2^{K-4}\vert \widehat C^{(3)}_{KH}\vert}{\vert C_{KH}\vert}\,\vert \widehat\rho_{KH}(j-k)\vert\,,
\end{align*}
for $\vert j-k\vert\gs 1$, where
\begin{align*}
C_{KH}=&(KH-1)KH(KH+1)(KH+2)\\
\widehat\rho_{KH}(j-k)=&6(j-k)^{KH+2}-4(j-k-1)^{KH+2}-4(j-k+1)^{KH+2}\\
&+(j-k-2)^{KH+2}+(j-k+2)^{KH+2}.
\end{align*}
So
\[
\big\vert \widehat d_{jk}^{B^{KH}}\big\vert\ls
2^{1-K}\vert\rho_{KH}(j-k)\vert \bigg(\frac{T}{n}\bigg)^{2KH} +\frac{2^{K-2}\vert \widehat C^{(2)}_{KH}\vert+2^{K-4}\vert \widehat C^{(3)}_{KH}\vert}{\vert C_{KH}\vert}\,\vert \widehat\rho_{KH}(j-k)\vert \bigg(\frac{T}{n}\bigg)^{2KH}.
\]
The cases when $j=1$ and $k>1$ or $k=1$ and $j>1$ can be proven in a way analogous to that of subfBm.
Next,  we obtain
\[
\max_{1\ls k\ls n-1}\sum_{j=1}^{n-1}\vert d^{B^{KH}}_{jk}\vert \ls 1+\Big(\max_{1\ls k\ls n-2}\vert d^{B^{KH}}_{k+1,k}\vert+\max_{2\ls k\ls n-1}\vert d^{B^{KH}}_{k-1,k}\vert\Big)
+\max_{1\ls k\ls n-1}\sum_{j=1\atop \vert j-k\vert\gs 2}^{n-1}\vert d^{B^{KH}}_{jk}\vert\,.
\]
Since
\[
\vert d_{k+1,k}^{B^{KH}}\vert\ls \sqrt{d_{k,k}^{B^{KH}}d_{k+1,k+1}^{B^{KH}}},
\]
then by the inequality (\ref{covariance1}) we get
\[
\max_{1\ls k\ls n-2}\vert d^{B^{KH}}_{k+1,k}\vert\ls 1.
\]
Reasoning as in \cite{bc}, Appendix A we obtain
\[
\sum_{j=2}^\infty\vert\widehat\rho_{KH}(j)\vert=3-3\cdot 2^{KH+2}+3^{KH+2}.
\]
Furthermore,
\begin{align*}
&\max_{1\ls k\ls n-1}\sum_{j=1\atop \vert j-k\vert\gs 2}^{n-1}\vert d^{B^{KH}}_{jk}\vert\\
&\quad\ls \max_{1\ls k\ls n-1}\sum_{j=1\atop \vert j-k\vert\gs 2}^{n-1}\frac{\vert\rho_{KH}(j-k)\vert}{4-4^{KH}}+ \frac{2^{2K-3}\vert \widehat C^{(2)}_{KH}\vert+2^{2K-5}\vert \widehat C^{(3)}_{KH}\vert}{\vert C_{KH}\vert}\,\max_{1\ls k\ls n-1}\sum_{j=1\atop \vert j-k\vert\gs 2}^{n-1}\frac{\vert \widehat\rho_{KH}(j-k)\vert}{4-4^{KH}}\\
&\quad\ls -\frac 12\,\sgn(2KH-1)\,\frac{3-3\cdot 2^{2KH}+3^{2KH}}{4-4^{KH}}+\frac{2^{2K-3}\vert \widehat C^{(2)}_{KH}\vert+2^{2K-5}\vert \widehat C^{(3)}_{KH}\vert}{\vert C_{KH}\vert}\,\,\frac{3-3\cdot 2^{KH+2}+3^{KH+2}}{4-4^{KH}}\,.
\end{align*}
It is clear that
\begin{align*}
\frac{2^{2K-3}\vert \widehat C^{(2)}_{KH}\vert}{\vert C_{KH}\vert}=&\frac{(1-K)(2-K)(3-K)2^4 H^3 2^{K-3}}{(1-KH) (KH+1)(KH+2)}=\frac{(H-HK)(2-K)(3-K)2^4 H^2 2^{K-3}}{(1-KH) (KH+1)(KH+2)}\\
\ls&\frac{(2-K)(3-K) H^2 2^{K+1}}{ (KH+1)(KH+2)}\ls
\frac{6\cdot 2^{K-1}}{2}\ls 3
\end{align*}
and
\begin{align*}
\frac{2^{2K-5}\vert \widehat C^{(3)}_{KH}\vert}{\vert C_{KH}\vert}=&\frac{2^{K-3}(1-K) H \vert 2H-1\vert (1+2H -2KH)}{(1-KH)(KH+1)(KH+2)}\\
\ls& \frac{2^{K-3} \vert 2H-1\vert (1+2H -2KH)}{(KH+1)(KH+2)}\ls
\frac{2\cdot 2^{K-3}}{2}\ls \frac 14\,.
\end{align*}
Thus
\begin{align*}
\max_{1\ls k\ls n-1}\sum_{j=1}^{n-1}\vert d^{B^{KH}}_{jk}\vert
\ls& 3+\frac{3-3\cdot 2^{2KH}+3^{2KH}}{2(4-4^{KH})}+\frac{13}{4}\cdot\frac{3-3\cdot 2^{KH+2}+3^{KH+2}}{4-4^{KH}}\\
\ls& 3+\frac 14+\frac{13}{4}\cdot\frac{54}{100}=5,005
\end{align*}
since the numerator of the second term is a decreasing function, and the numerator of the third term is an increasing function.

\subsection{Ornstein-Uhlenbeck process}

The \textbf{fractional Ornstein-Uhlenbeck} (fO-U) process of the first kind is the unique solution of the stochastic differential equation
\begin{equation}\label{O-U1}
X_t=x_0-\mu\int_0^t X_s\,ds+\t B^H_t,\qquad t\ls T,
\end{equation}
with $\mu,\t>0$, where $B^H$, $0<H<1$, is a fBm. Its explicit solution is given by
\[
X_t=x_0 e^{-\mu t}+\t\int_0^t e^{-\mu(t-u)}dB^H_u,
\]
where the integral exists as a Riemann-Stieltjes integral for all $t > 0$ (see, e.g., \cite{chm}).
First we show the following lemma.
\begin{lem}\label{lem}
Let $X$ be the solution of equation (\ref{O-U1}). Assume that $B^H=\{B^H(t): t\in[0,T]\}$ is a fBm with  $H\in(0,H^*]$, where a real number  $H^*<1$ is known. Then
for $H\neq 1/2$
\[
\max_{1\ls k\ls n-1}\sum_{j=1}^{n-1}  \vert  d^{X,n}_{jk}\vert\ls\frac {4\mu T}{4-2^{2H^*}}\bigg( 2\mu\,\frac{T}{n}\,\bigg[C \bigg(\frac Tn\bigg)^{2-2H^*} +1\bigg]+\sqrt{2(4-2^{2H^*})}\,\bigg[C \bigg(\frac Tn\bigg)^{2-2H^*} +1\bigg]^{1/2}\bigg)+\frac 83\,,
\]
where $C=\mu^2\big[3x_0^2\t^{-2}+6T^2\big]$ and
\[
d^{X,n}_{jk}=\frac {n^{2H}}{\t^2 T^{2H}(4-2^{2H})}\,\widehat d^{X,n}_{jk},\qquad \widehat d^{X,n}_{jk}=\E (\D_{n,j}^{(2)} X \D_{n,k}^{(2)} X).
\]
\end{lem}
\proof Denote
\[
\D_{n,k} Z=-\mu\int_{t^n_{k-1}}^{t^n_k} X_s\,ds.
\]
It is clear that
\begin{align*}
\big\vert \widehat d^{X,n}_{jk}-\t^2  \widehat d^{B^H,n}_{jk}\big\vert
=& \big\vert\E [\D_{n,j+1} Z-\D_{n,j} Z][\D_{n,k+1} Z-\D_{n,k} Z] +\t\E [\D_{n,j+1} Z-\D_{n,j} Z]\D^{(2)}_{n,k}B^H\\
&+\t\E [\D_{n,k+1} Z-\D_{n,k} Z]\D^{(2)}_{n,j}B^H\big\vert
\end{align*}
and
\[
[\D_{n,j+1} Z-\D_{n,j} Z]=-\mu\bigg[\int_{t^n_j}^{t^n_{j+1}} [X_s-X_j]\,ds-\int_{t^n_{j-1}}^{t^n_j} [X_s-X_j]\,ds\bigg].
\]
Reasoning as in \cite{kk1}  we obtain
\[
\sup_{t\ls T}\E X_t^2\ls 3x_0^2+6\t^2 T^2
\]
and
\begin{align*}
\E [X_t-X_j]^2\ls&2\mu^2\frac Tn \int_{t^n_j}^t \E X^2_s\,ds+2\t^2 \bigg(\frac Tn\bigg)^{2H}\ls 2\mu^2\big[3x_0^2+6\t^2 T^2\big]\bigg(\frac Tn\bigg)^2 +2\t^2 \bigg(\frac Tn\bigg)^{2H}\\
=&2\t^2\bigg(\frac Tn\bigg)^{2H}\bigg[C \bigg(\frac Tn\bigg)^{2-2H} +1\bigg],
\end{align*}
where $C=\mu^2\big[3x_0^2\t^{-2}+6T^2\big]$.
Thus
\begin{align*}
\big\vert \widehat d^{X,n}_{jk}-\t^2 \widehat d^{B^H,n}_{jk}\big\vert
\ls&8\mu^2\t^2 \bigg(\frac Tn\bigg)^{2+2H}\bigg[C \bigg(\frac Tn\bigg)^{2-2H} +1\bigg]\\
&+ 4\sqrt{2(4-2^{2H})}\,\mu\t^2\bigg(\frac Tn\bigg)^{1+2H}\bigg[C \bigg(\frac Tn\bigg)^{2-2H} +1\bigg]^{1/2}.
\end{align*}
Consequently,
\begin{align*}
&\max_{1\ls k\ls n-1}\sum_{j=1}^{n-1}  \vert d^{X,n}_{jk}\vert\ls \max_{1\ls k\ls n-1}\sum_{j=1}^{n-1}\big\vert d^{X,n}_{jk}-\t^2 d^{B^H,n}_{jk}\big\vert+\max_{1\ls k\ls n-1}\sum_{j=1}^{n-1} \big\vert  d^{B^H,n}_{jk}\big\vert\\
\ls& \frac {8\mu^2}{4-2^{2H^*}}\frac{T^2}{n}\bigg[C \bigg(\frac Tn\bigg)^{2-2H^*} +1\bigg]+\frac {4\sqrt{2}\,\mu T}{\sqrt{4-2^{2H^*}}}\,\bigg[C \bigg(\frac Tn\bigg)^{2-2H^*} +1\bigg]^{1/2}+\frac 83\\
=&\frac {4\mu T}{4-2^{2H^*}}\bigg( 2\mu\,\frac{T}{n}\,\bigg[C \bigg(\frac Tn\bigg)^{2-2H^*} +1\bigg]+\sqrt{2(4-2^{2H^*})}\,\bigg[C \bigg(\frac Tn\bigg)^{2-2H^*} +1\bigg]^{1/2}\bigg)+\frac 83\,.
\end{align*}

\begin{lem} Let $X$ be the solution of the equation (\ref{O-U1}). Then
\begin{align*}
\bigg\vert\frac{1}{n-1}\sum_{k=1}^{n-1} d^{X,n}_{kk}-1\bigg\vert\ls& \frac {4\mu}{4-2^{2H^*}}\,\frac {T}{n}\bigg(2\mu\,\frac Tn\,\bigg[C \bigg(\frac Tn\bigg)^{2-2H^*} +1\bigg]\\
&+\sqrt{2(4-2^{2H^*})}\,\bigg[C \bigg(\frac Tn\bigg)^{2-2H^*} +1\bigg]^{1/2}\bigg).
\end{align*}
\end{lem}
\proof To prove this lemma, observe that
\begin{align*}
\bigg\vert\frac{1}{n-1}\sum_{k=1}^{n-1} d^{X,n}_{kk}-1\bigg\vert\ls& \max_{1\ls k\ls n-1}\big\vert d^{X,n}_{kk}- d^{B^{H,n}}_{kk}\big\vert +\bigg\vert\frac{1}{n-1}\sum_{k=1}^{n-1} d^{B^{H,n}}_{kk}-1\bigg\vert\\
=& \max_{1\ls k\ls n-1}\big\vert d^{X,n}_{kk}-1\big\vert\,.
\end{align*}
From the inequality
\begin{align*}
\big\vert \widehat d^{X,n}_{kk}-\t^2 \widehat d^{B^H,n}\big\vert
\ls&8\mu^2\t^2 \bigg(\frac Tn\bigg)^{2+2H}\bigg[C \bigg(\frac Tn\bigg)^{2-2H} +1\bigg]\\
&+ 4\sqrt{2(4-2^{2H})}\,\mu\t^2\bigg(\frac Tn\bigg)^{1+2H}\bigg[C \bigg(\frac Tn\bigg)^{2-2H} +1\bigg]^{1/2}
\end{align*}
it follows that the statement of the lemma holds.

\section{Simulations}\label{s:modeling}
\setcounter{figure}{0}

The simulations of the obtained confidence intervals presented below were performed using the R software environment \cite{R}. Sample paths of fBm were generated using the circulant matrix embedding method and were further used to simulate the sample paths of the fractional Ornstein-Uhlenbeck process (\ref{O-U1}). The constants for the latter were (arbitrarily) chosen as $x_0 = 0$ and $\mu=0.5$. Sample paths of the sub-fractional and bifractional Brownian motion were simulated using the Cholesky method. Due to the notable computational requirements of this method the maximum sample path length considered was $n=1600$. Figures presented below correspond to the case of the confidence level $1-\alpha$, $\alpha=0.1$. The observed coverage percentages in all cases were at least as good as claimed in Theorem \ref{pagrindine}.

\begin{figure}[h]
\centering
\begin{minipage}{0.45\textwidth}
\centering
\includegraphics[scale=0.6]{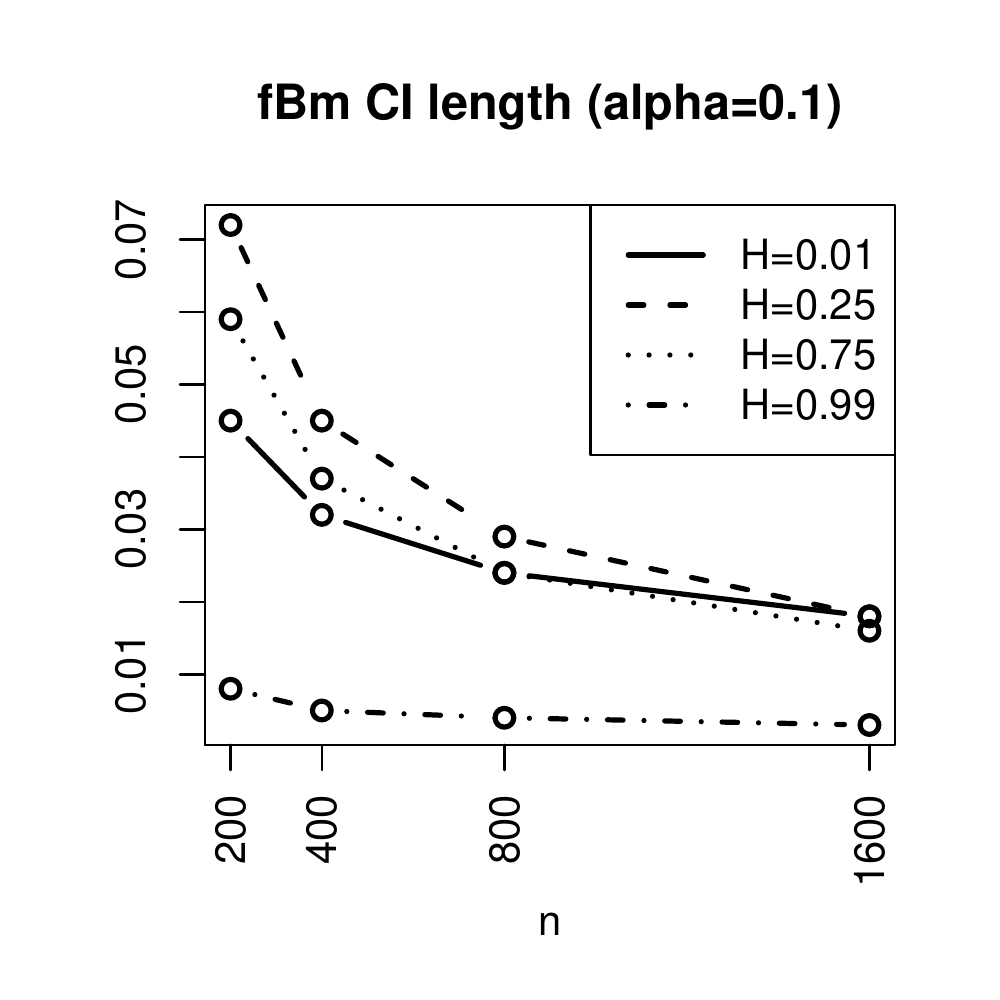}
\caption{$\big\vert H_n^{sup}(\alpha) - H_n^{inf}(\alpha) \big\vert $}
\label{fig:c_fbm}
\end{minipage}\hfill
\begin{minipage}{0.45\textwidth}
\centering
\includegraphics[scale=0.6]{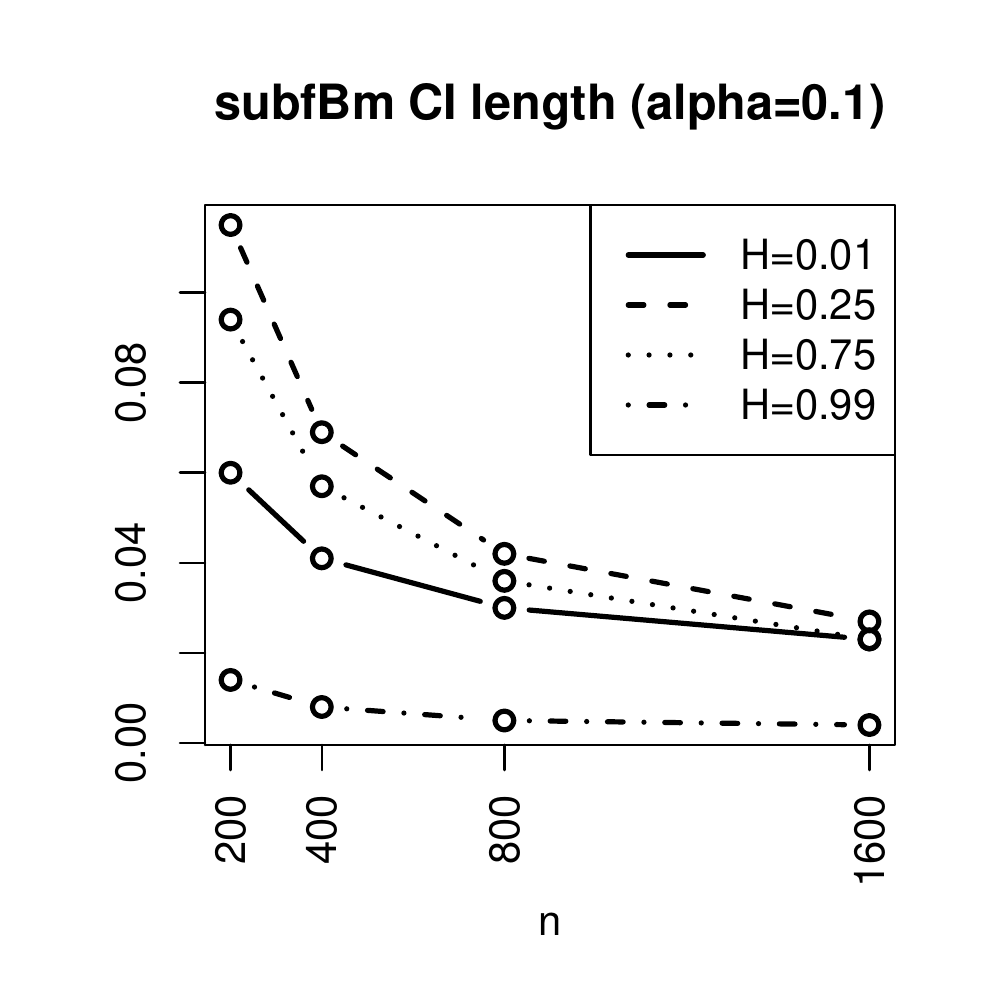}
\caption{$\big\vert H_n^{sup}(\alpha) - H_n^{inf}(\alpha) \big\vert $}
\label{fig:c_subfbm}
\end{minipage}
\end{figure}

\begin{figure}[h]
\centering
\begin{minipage}{0.45\textwidth}
\centering
\includegraphics[scale=0.6]{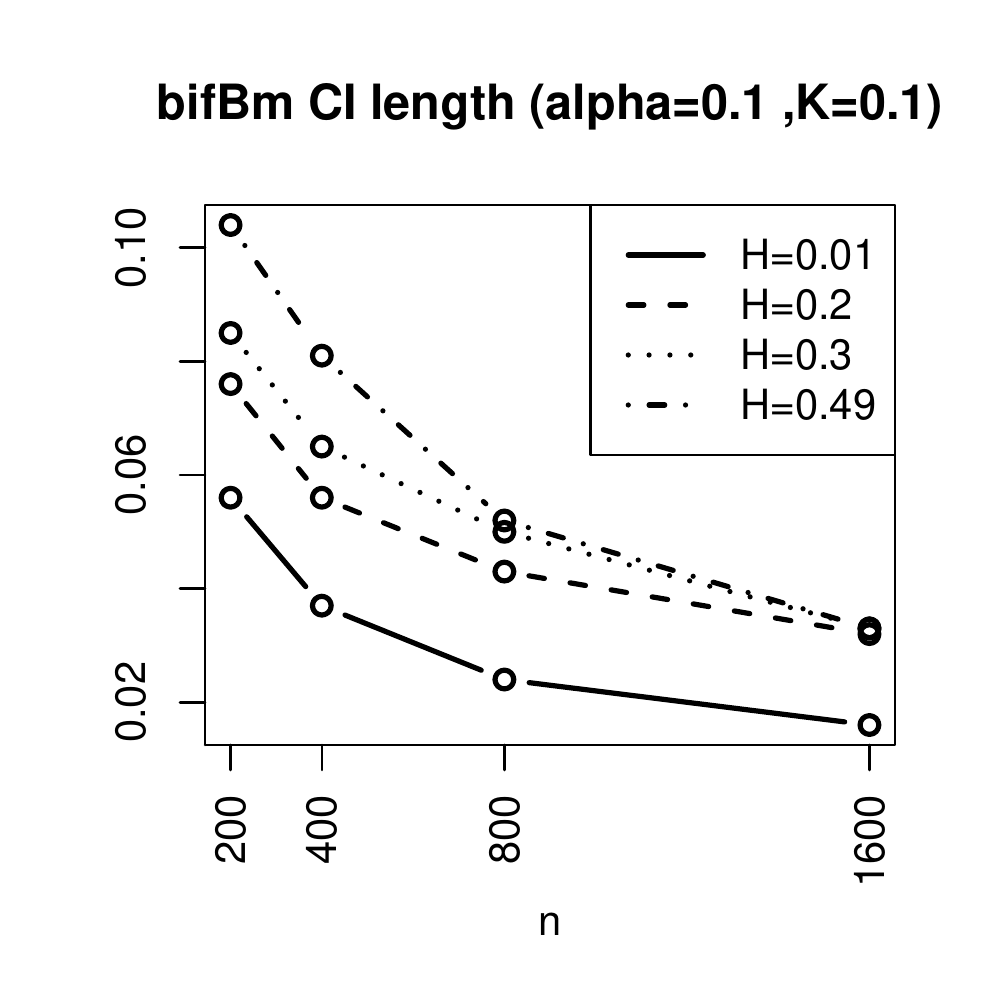}
\caption{$\big\vert (HK)_n^{sup}(\alpha) - (HK)_n^{inf}(\alpha) \big\vert $}
\label{fig:c_bifbm_lowK}
\end{minipage}\hfill
\begin{minipage}{0.45\textwidth}
\centering
\includegraphics[scale=0.6]{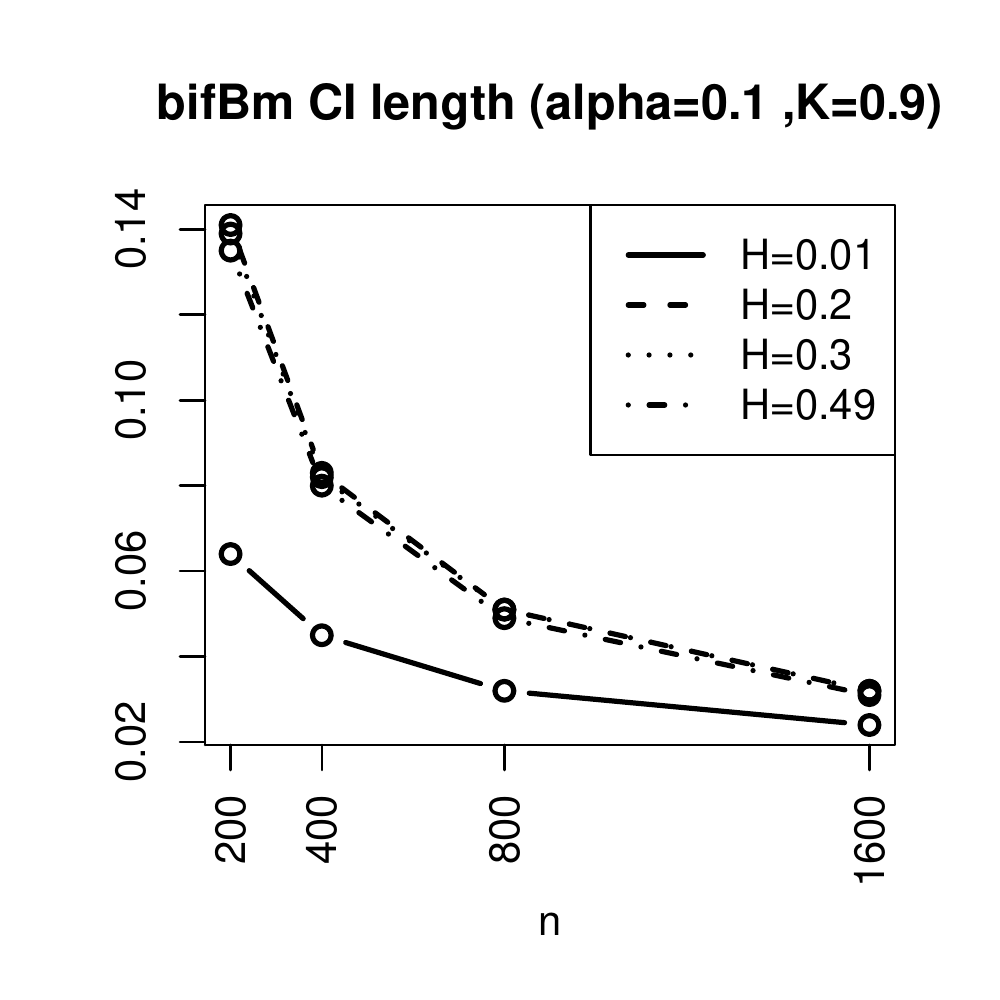}
\caption{$\big\vert (HK)_n^{sup}(\alpha) - (HK)_n^{inf}(\alpha) \big\vert $}
\label{fig:c_bifbm_highK}
\end{minipage}
\end{figure}

\begin{figure}[h]
\centering
\begin{minipage}{0.45\textwidth}
\centering
\includegraphics[scale=0.6]{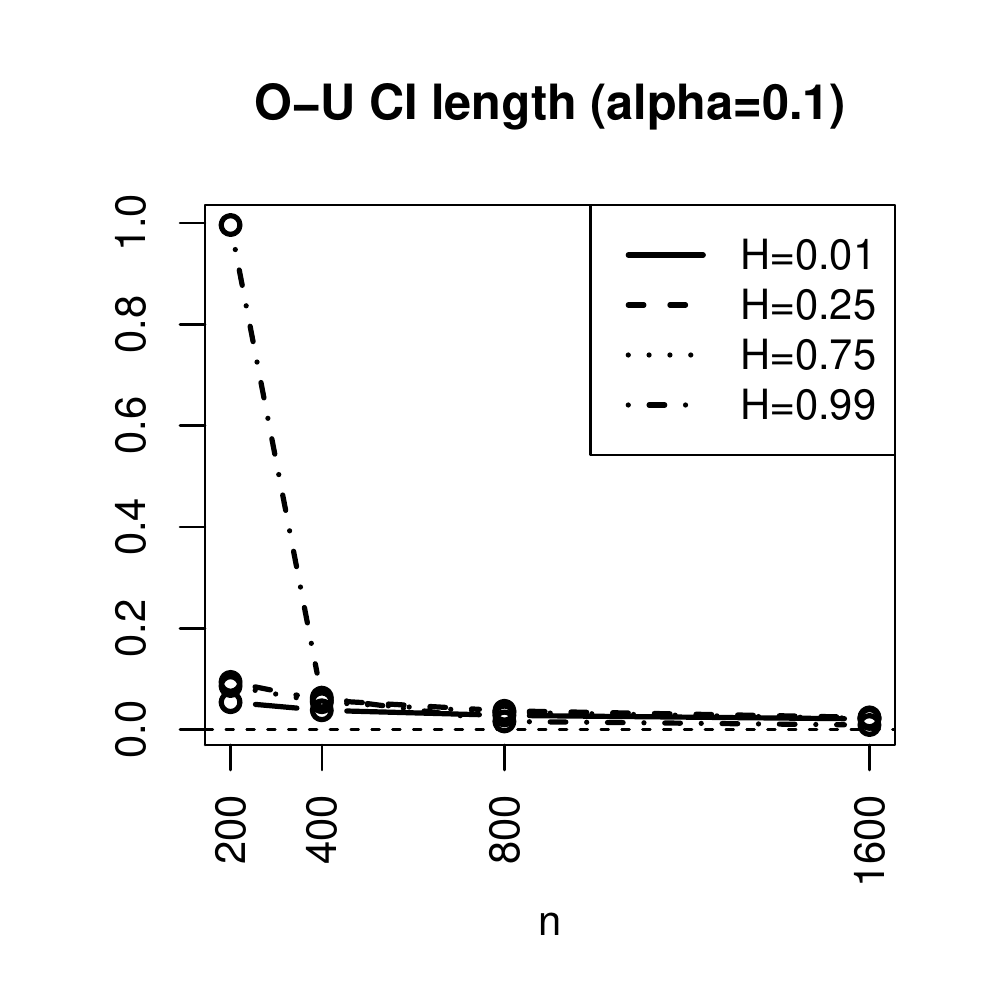}
\caption{$\big\vert H_n^{sup}(\alpha) - H_n^{inf}(\alpha) \big\vert $}
\label{fig:c_OU}
\end{minipage}\hfill
\begin{minipage}{0.45\textwidth}
\centering
\includegraphics[scale=0.4]{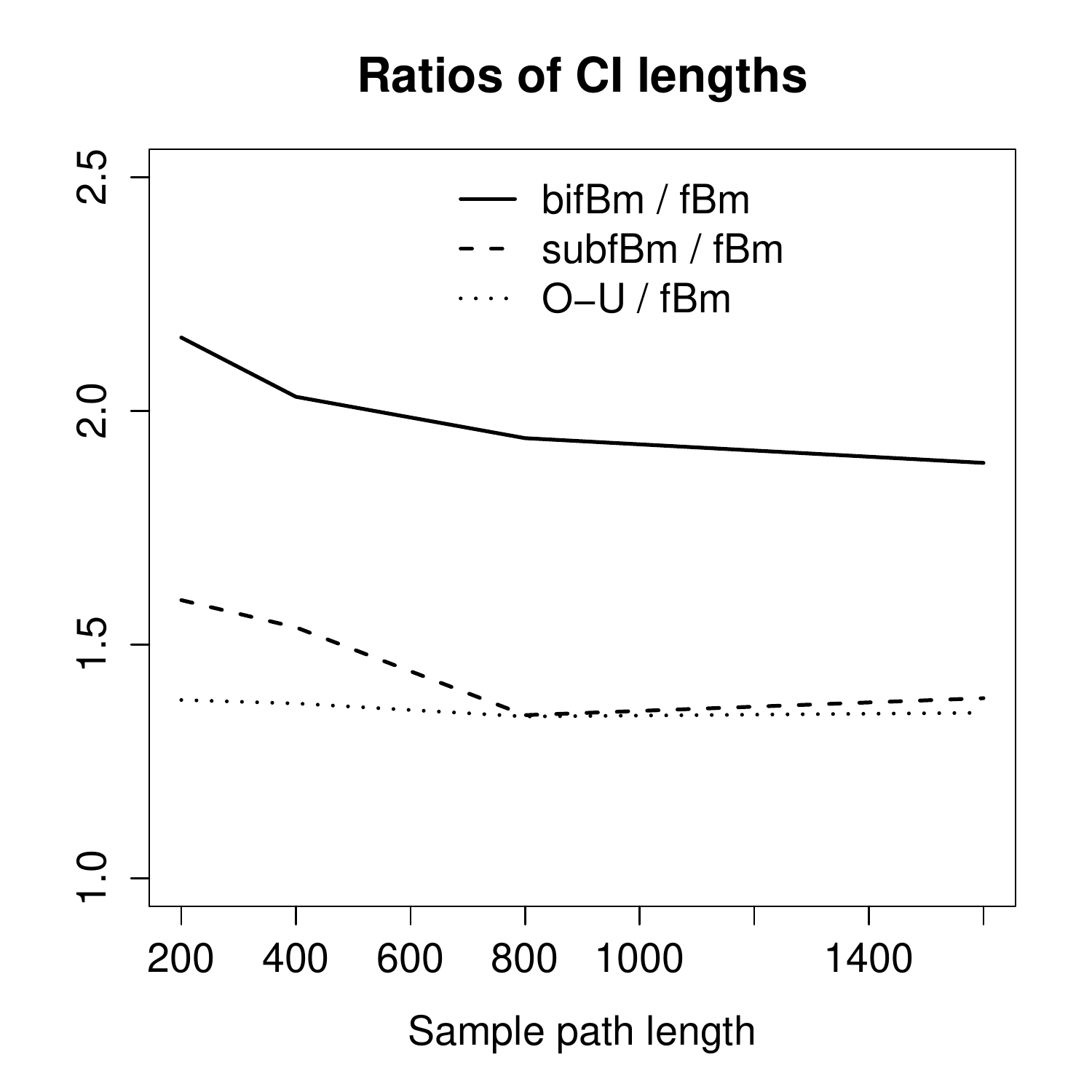}
\caption{Comparison of CI lengths}
\label{fig:c_length_ratios}
\end{minipage}
\end{figure}

Figures \ref{fig:c_fbm} - \ref{fig:c_OU} present the confidence interval (CI) lengths for all the process types considered in this paper. Figure \ref{fig:c_length_ratios} shows the median ratios of the confidence intervals lengths, where the CI lengths of the subfBm, bifBm and fO-U processes were divided by the corresponding CI lengths of fBm. It can be seen that in almost all cases the confidence intervals behave in a similar way, one notable exception being the case of fO-U as the value of $H$ approaches 1. This is hardly unexpected given the normalization used in Lemma \ref{lem}, and in this scenario the CI covers the whole interval of possible parameter values $0 < H < 1$.

\section*{Appendix. Code listings}

\subsection*{genFBM.r}

\begin{lstlisting}[language=R]
##################################################
# genFBM returns a single sample path of the fBm #
# H - the Hurst index                            #
# N - length of the sample path                  #
##################################################

genFBM <- function(H,N) {
  # M - length of the required vector of autocovariances
  M <- 2^(ceiling( 1 + log(N-1)/log(2) ))
  # V is the M-vector of autocovariances of fGn
  H2 <- 2*H
  V <- c(0:(M/2), (M/2-1):1)
  V <- 1/(2*N^(H2))*(abs(V-1)^(H2) - 2*V^(H2) + (V+1)^(H2))
  # W is the fast Fourier transform of V
  W <- Re(fft(V))
  # We increase M until all coordinates of W are positive
  while (any(W<=0) & (M<2^30)) {
    M <- 2*M
    V <- c(0:(M/2), (M/2-1):1)
    V <- 1/(2*N^(H2))*(abs(V-1)^(H2) - 2*V^(H2) + (V+1)^(H2))
    W <- Re(fft(V))
  }
  # X and Y and iid Gaussian with mean 0 and variance 1/sqrt(2)
  X <- rnorm(M, mean=0, sd=(1/sqrt(2)))
  Y <- rnorm(M, mean=0, sd=(1/sqrt(2)))
  Z <- vector(length=M)
  Z[1] <- X[1]
  Z[M/2+1] <- Y[1]
  Z[c(2:(M/2))] <- X[c(2:(M/2))] + 1i*Y[c(2:(M/2))]
  Z[M+2-c(2:(M/2))] <- X[c(2:(M/2))] - 1i*Y[c(2:(M/2))]
  U <- sqrt(W)*Z
  # X is the fast Fourier transform of U = sqrt(W) * ( X + iY )
  X <- fft(U, inverse=T)
  # The real part of the first N coordinates of X is the desired
  # sample path, in this case, the fGn.
  # We calculate the fBm sample path as cumulated sums of fGn.
  BHinc <- c(0,Re(X[1:N]))/sqrt(M)
  BH <- cumsum(BHinc)
  return(BH[-length(BH)])
}
\end{lstlisting}

\subsection*{Cholesky.r}

\begin{lstlisting}[language=R]
###############################################
# genSFBMch generates a batch of sample paths #
#   of sub-fractional Brownian motion         #
# H - the Hurst index                         #
# N - length of sample paths                  #
# Nsp - number of sample paths                #
###############################################

genSFBMch <- function(H,N,Nsp) {
  # sfBm's covariation matrix
  sfBm_cov <- function(t,s,H) {
    H2 <- H*2
    cts <- s^(H2)+t^(H2)-((t+s)^(H2)+abs(t-s)^(H2))/2
    return(cts)
  }

  G <- matrix(ncol=N,nrow=N)
  for (j in c(1:N)) {
    i <- seq(1:N)/N
    G[,j] <- (apply(as.matrix(i),1,sfBm_cov,j/N,H))
    }
  # Cholesky decomposition
  L <- t(chol(G))

  # a function to generate a single sample path
  gen_sp <- function(L) {
    Z <- rnorm(N)
    B <- c(0,(L %*% Z))
    return(B)
  }

  # generating Nsp of them
  BHM <- replicate(Nsp,gen_sp(L))
  return(BHM)
}

###############################################
# genSFBMch generates a batch of sample paths #
#   of bifractional Brownian motion           #
# H, K - parameters of the bifBm              #
# N - length of sample paths                  #
# Nsp - number of sample paths                #
###############################################

genBFBMch <- function(H,K,N,Nsp) {
  # bifBm's covariation matrix
  bifBm_cov <- function(t,s,H) {
    H2 <- H*2
    cts <- ((abs(t)^(H2)+abs(s)^(H2))^K-abs(t-s)^(K*H2))/(2^K)
    return(cts)
  }

  G <- matrix(ncol=N,nrow=N)
  for (j in c(1:N)) {
    i <- seq(1:N)/N
    G[,j] <- (apply(as.matrix(i),1,bifBm_cov,j/N,H))
    }
  # Cholesky decomposition
  L <- t(chol(G))

  # a function to generate a single sample path
  gen_sp <- function(L) {
    Z <- rnorm(N)
    B <- c(0,(L %*% Z))
    return(B)
  }

  # generating Nsp of them
  BHM <- replicate(Nsp,gen_sp(L))
  return(BHM)
}
\end{lstlisting}

\subsection*{fbmCI.r}

The following code evaluates the confidence intervals for the fractional Brownian motion. Simulations for other processes considered in this paper were performed in a similar way.

\begin{lstlisting}[language=R]
source("Cholesky.r")
source("genFBM.r")

############################
# Definitions of functions #
############################

phiL = function(z, nu, eps, N) {
  if (z >= 0 & z < (eps+1) * sqrt(N-1)) {
    return ( exp(z*sqrt(N-1) / (2*nu))
             * (1 - z/((eps+1) * sqrt(N-1)))
             **( (eps+1) * (N-1) / (2*nu) ) ) }
  else return(0)
}

phiR = function(z, nu, eps, N) {
  return ( exp(-z*sqrt(N-1)/ (2*nu))
           * (1 + z/((eps+1) * sqrt(N-1)))
           **( (eps+1) * (N-1) / (2*nu) ) )
}

phiLOpt = function(z, nu, eps, N, value){
  return(abs(phiL(z, nu, eps, N) - value))
}

phiROpt = function(z, nu, eps, N, value){
  return(abs(phiR(z, nu, eps, N) - value))
}

gn = function(x, N) {
  if (x == 0) { return(-Inf) } else { return(2*x*log(N) - log(4-4**x)) }
}

gnOpt = function(x, N, value) {
  return(abs(gn(x, N) - value))
}

logStar = function(x) {
  if (x>0) { return(log(x)) } else { return(-Inf) }
}

confBC = function(B, alpha, nu, eps, N, kappa) {
  qL = optimize(f=phiLOpt, interval=c(0, 20), value=alpha/2,
                nu=nu, eps=eps, N=(N-1))[[1]]
  qR = optimize(f=phiROpt, interval=c(0, 20), value=alpha/2,
                nu=nu, eps=eps, N=(N-1))[[1]]

  xL = 1-qL/sqrt(N-1)
  xR = 1+qR/sqrt(N-1)
  Sn = mean(diff(B, differences=2)**2)

  argL = max( logStar((xL - eps)*kappa**2 / Sn), -log(3) )
  argR = log((xR + eps)*kappa**2 / Sn)

  Hmin = max(c(0, optimize(gnOpt, c(0.00001, 0.99999),
             N=N, value=argL)[[1]]))
  Hmax = min(c(1, optimize(gnOpt, c(0.00001, 0.99999),
             N=N, value=argR)[[1]]))

  return(c(Hmin, Hmax))
  }

#################################
# Parameters used for modelling #
#################################

nseq = c(200, 400, 800, 1600)
Hseq = c(0.01, 0.25, 0.75, 0.99)
Nsp = 1000
alpha = 0.1
kappa = 1

############################################
# ciLenM - matrix of CI lengths            #
# ciInsideM - matrix of CI coverage ratios #
############################################

ciLenM = matrix(nrow=length(Hseq), ncol=length(nseq))
ciInsideM = matrix(nrow=length(Hseq), ncol=length(nseq))
for (m in 1:length(Hseq)) {
  H = Hseq[m]
  ciLen <- vector(length=length(nseq))
  ciInside <- vector(length=length(nseq))
  for (k in 1:length(nseq)) {
    N = nseq[k]
    BM = replicate(Nsp, genFBM(H, N))
    eps = 1/N
    nu = 8/3
    CI = apply(BM, 2, confBC, alpha=alpha, nu=nu, eps=eps,
               N=N, kappa=kappa)
    ciLen[k] <- mean(CI[2,] - CI[1,])
    ciInside[k] <- length(which(CI[1,]<H & CI[2,]>H))/Nsp
  }
  ciLenM[m,] = round(ciLen, 3)
  ciInsideM[m,] = ciInside
}

########################
# Plotting the results #
########################

pdf("bw_fbm_ci.pdf", width = 8, height = 4)
  par(mfrow=c(1,2))
  yl = c(min(ciLenM), max(ciLenM))
  plot(nseq, ciLenM[1,], type="b", lwd=2, lty=1, xlab="n", ylab="",
       main=paste0("fBm CI length (alpha=", alpha,")"), ylim=yl, xaxt="n")
  for (k in 2:length(Hseq)) {
    lines(nseq, ciLenM[k,], type="b", lwd=2, lty=k)
  }
  abline(h=0, lty=2)
  legend("topright", c("H=0.01", "H=0.25", "H=0.75", "H=0.99"),
         lwd=2, lty=1:length(Hseq))
  axis(1, nseq, labels=nseq, las=2)

  yl = c(min(ciInsideM), max(ciInsideM))
  plot(nseq, ciInsideM[1,], type="b", lty=1, xlab="n", ylab="",
       main="Coverage %", ylim=yl,  xaxt="n")
  for (k in 2:length(Hseq)) {
    lines(nseq, ciInsideM[k,], type="b", lty=k)
  }
  legend("bottomright",c("H=0.01", "H=0.25", "H=0.75", "H=0.99"),
         lwd=2, lty=1:length(Hseq))
  axis(1, nseq, labels=nseq, las=2)
dev.off()
\end{lstlisting}

\end{document}